\documentclass[12pt]{article}  % `article' with A4 paper size option
\usepackage{amsmath,amssymb}

\newtheorem{theorem}{Theorem}

\newenvironment{remark}
{\smallskip\noindent{\bf Remark\/}.}{\smallskip\par}
\newenvironment{remarks}
{\smallskip\noindent{\bf Remarks\/}.}{\smallskip\par}
\newenvironment{proof}
{\noindent{\em Proof\/}.}{{ $\Box$}\smallskip\par}

\newcommand{\CC}{{\Bbb C}}

\newcommand{\PP}{{\Bbb P}}

\newcommand{\RR}{{\Bbb R}}

\newcommand{\VV}{{\Bbb V}}

\newcommand{\calM}{{\cal M}}

\newcommand{\eps}{\varepsilon}

\title{Deformations of polynomials and their zeta functions
}
\author{S.~M.~Gusein-Zade 
\thanks{Partially supported by the grants RFBR--04--01--00762,
NSh--1972.2003.1}
\and D.~Siersma
}

\date{}

\begin{document}

\maketitle

\begin{abstract}
For an analytic in $\sigma\in(\CC, 0)$ family $P_\sigma$ of polynomials in $n$ variables
there is defined a monodromy transformation $h$ of the zero level set $V_\sigma=\{P_\sigma=0\}$
for $\sigma\ne 0$ small enough. The zeta function of this monodromy transformation
is written as an integral with respect to the Euler characteristic of the corresponding
local data. This leads to a study of deformations of holomorphic germs and their zeta functions.
We show some examples of computations with the use of this technique.
\end{abstract}

\section{Introduction}
A complex polynomial $P$ in $n$ variables defines a map from $\CC^n$ to $\CC$
(also denoted by $P$). This map is a ($C^\infty$) locally trivial fibration over 
the complement to a finite subset of the target $\CC$~--- the bifurcation set:
\cite{Varchenko}. Topological study of polynomial maps was started in
\cite{Broughton} and continued, in particular, in \cite{Parus, ST1}. One is
interested in a description of the bifurcation set of a polynomial, its generic
level set, degeneration of the level set for singular values, monodromies 
around
singular values, {\dots} The level sets $\{P=\sigma\}$ of the polynomial $P$ are
the zero level sets of the family of polynomials $P_\sigma=P-\sigma$. It is
natural to study the behaviour of the zero level sets not of this particular
family, but of a general family $P_\sigma$ of polynomials (say, analytic in
$\sigma$ from a neighbourhood of zero in the complex line $\CC_\sigma$). In this
setting one can also study changes of the map $P_\sigma: \CC^n \to \CC$ in the
family. Such a study was started in \cite{ST2, ST3}.

Let $P_\sigma(x)$ be an analytic in $\sigma\in(\CC_\sigma, 0)$
family of polynomials in $n$ variables $x=(x_1, \ldots, x_n)$:
a deformation of the polynomial $P=P_0$. The polynomial $P_\sigma$
defines a map $P_\sigma:\CC^n \to \CC$. Let $\VV$ be the hypersurface
in $\CC^{n+1}=\CC^n\times\CC_\sigma$ defined by the equation
$P_\sigma(x)=0$ ($x\in\CC^n$, $\sigma\in\CC_\sigma$). The projection
$p:\VV\to\CC_\sigma$ to the second factor is a fibre bundle over
a punctured neighbourhood of the origin in $\CC_\sigma$ (see, e.g.,
\cite{Varchenko}). Let $V_\sigma=\{P_\sigma=0\} \subset\CC^n$, $\sigma\ne0$
small enough, be its fibre (the zero level set of a generic member
of the family) and let $h:V_\sigma \to V_\sigma$ be a monodromy
transformation of the described fibration (it is well defined up to
isotopy).

For a transformation $h:Z\to Z$ of a space $Z$, its zeta function
$\zeta_h(t)$ is the rational function
$$
\prod_{q\ge 0} \{\det(id-t\cdot h_{*\vert H_q(Z;\RR)})\}^{(-1)^q}.
$$
The degree of the zeta function (the degree of the numerator
minus the degree of the denominator) is equal to the Euler
characteristic $\chi(Z)$ of the space $Z$.
The zeta function of the monodromy transformation $h:V_\sigma \to V_\sigma$
will be called the zeta function of the family $P_\sigma$.

For $P_\sigma=P+\sigma$, $V_\sigma$ is a generic level set of the polynomial $P$
and $\zeta_P(t)=\zeta_h(t)$ is the zeta function of the classical monodromy
transformation of the polynomial function $P$ around the zero value. In
\cite{GLM3}, there is a formula which expresses the zeta function $\zeta_P(t)$
(and thus the Euler characteristic of a generic level set of the polynomial $P$)
as an integral with respect to the Euler characteristic of the corresponding
local data at each point of the compactification $\CC\PP^n$ of the affine space
$\CC^n$. This "localization" appeared to be effective for computing these
invariants in a number of cases.

Here we formulate a general localization formula for the zeta function and
specialise it for families of polynomials. For that we describe local data
corresponding to the problem of computing the zeta function of a family of
polynomials. These are deformations of germs of (complex analytic) functions on
the affine space $\CC^n$ or on the affine space $\CC^n$ with a distinguished
hyperplane ("boundary") $\CC^{n-1}$ and their zeta functions. We show some
examples of application of the localization formula.

\section{The localization principle}
Let $X$ be a compact complex analytic (generally speaking
singular) variety and let $Y$ be a (compact) subvariety of $X$. Let
$L$ be a line bundle over $X$ and let $s_\sigma$ be an
analytic in $\sigma\in(\CC_\sigma,0)$ family of sections of the line
bundle $L$. Let $\VV\subset X\times \CC_\sigma$ be defined by the
equation $s_\sigma(x)=0$. The restriction of the projection
$p:X\times\CC_\sigma\to\CC_\sigma$ to the second factor to the
complement $\VV\setminus(Y\times\CC_\sigma)$ is a locally trivial
fibration over a punctured neighbourhood of the origin in
$\CC_\sigma$: \cite{Varchenko}. Let $V_\sigma=p^{-1}(\sigma)
\cap \left(\VV\setminus(Y\times\CC_\sigma)\right)$ ($\sigma\ne 0$
small enough) be its fibre (the zero set in $X\setminus Y$ of a
generic section from the family) and let $h:V_\sigma\to
V_\sigma$ be a monodromy transformation of the fibration (it is
well defined up to isotopy). Let $\zeta_{s_\sigma}(t)$ be the zeta function of
the monodromy transformation $h$ of the described fibration.

Now let's describe local versions of these objects. Pay attention that over a
neighbourhood of a point a line bundle is trivial and therefore its sections
can be considered as functions. Let $(X, x_0)\subset(\CC^N, x_0)$ be the germ of a 
complex analytic
variety and let $(Y, x_0)$ be a subvariety of it (possibly the empty one). Let
$s_\sigma$ be an analytic in $\sigma\in(\CC_\sigma,0)$ deformation of a germ $s$ 
of a function on $(X, x_0)$, i.e. $s_\sigma=S(\cdot, \sigma)$ where $S$ is a 
germ of a holomorphic function on $(X \times\ \CC_\sigma, (x_0,0)\,)$, $s_0=s$.
Let $\VV_{x_0}\subset (X\times \CC_\sigma, (x_0,0)\,)$ be the germ of the 
variety
defined by the equation $s_\sigma(x)=0$. Let positive $\eps$ be small enough
so that all the strata of a Whitney stratification of the pair $(X, \{s=0\})$
are transversal to the sphere $S_{\eps'}(x_0)$ of radius $\eps'$ centred
at the point $x_0$ in $\CC^N$ for any positive $\eps'\le\eps$.
Let $B_\eps(x_0)$ be the ball of radius $\eps$ centred at the point $x_0$.
The restriction of the projection $X\times\CC_\sigma\to \CC_\sigma$ to the
second factor to the complement $\left(\VV_{x_0}\cap(B_\eps\times\CC_\sigma)\right)
\setminus(Y\times\CC_\sigma)$ is a locally trivial
fibration over a punctured neighbourhood of the origin in $\CC_\sigma$:
\cite{Le} ($B_\eps$ is the ball of radius $\eps$).
Let $V_{\sigma,x_0}$ and $\zeta_{s_\sigma}(t)$ be the fibre (the
local zero set in $(X\setminus Y,0)$ of a generic function from the family) and
the zeta function of a monodromy transformation of this fibration.

For a constructible function $\Psi$ on a constructible set $Z$ with values in an
Abelian group $A$, there exists a notion of the integral $\int_Z\Psi d\chi$ of
the function $\Psi$ over the set $Z$ with respect to the Euler characteristic
(see, e.g., \cite{Viro}). For example, if $A$ is the multiplicative group of 
non-zero
rational functions in the variable $t$, $Z=\bigcup\Xi$ is a finite stratifiction
of $Z$ (without any regularity conditions) such that the function
$\Psi_x=\Psi(x)$ is one and the same for all points $x$ of each stratum $\Xi$ 
and
is equal to $\Psi_\Xi$ there, then by definition
$$
\int_Z\Psi d\chi = \prod_\Xi\left(\Psi_\Xi(t)\right)^{\chi(\Xi)}.
$$

\begin{theorem}\label{theo1}
One has
$$
\zeta_{s_\sigma}(t) = \int_{\{x\in X: s(x)=0\}}\zeta_{s_{\sigma, x}}(t)d\chi,
$$
and therefore
$$
\chi(V_\sigma) = \int_{\{x\in X: s(x)=0\}}\chi(V_{\sigma, x})d\chi.
$$
\end{theorem}

\begin{proof}
The monodromy transformation $h$ can be supposed to respect a Whitney
stratification of the pair $(X, Y)$. Because of the multiplicativity of the zeta
function, it is sufficient to proof the statement only for one stratum. Using
the induction one can suppose that the statement is already proved for strata of
lower dimension. Resolving singularities of the variety, we reduce the problem
to the case when $X$ is smooth. In this situation the proof is essentially the same
as in \cite{GLM1}.
\end{proof}

\begin{remark}
This statement can be considered as a generalization of the localization
principles described, in particular, in \cite{GLM1, GLM3} for particular cases.
It also can be deduced from general statements described in \cite{Dimca}.
\end{remark}

\section{Localization for families of polynomials}\label{pol}
Let $P_\sigma$ be an analytic in $\sigma\in(\CC_\sigma, 0)$ family of 
polynomials and
let $d$ be the degree of a generic polynomial of the family (i.e., the degree of
the polynomial $P_\sigma$ for $\sigma\ne 0$ small enough; $\deg P_0\le d$). Let
$X=\CC\PP^n$ be the (standard) compactification of the affine space $\CC^n$ and 
let
$Y=\CC\PP^{n-1}_\infty$ be its infinite hyperplane. The family of polynomials
$P_\sigma$ can be considered as a family of sections of the line bundle ${\cal
O}(-d)$ over the projective space $\CC\PP^n$ (if the degree of the polynomial 
$P_0$ is smaller than $d$, then the corresponding section vanishes on the whole
infinite hyperplane $\CC\PP^{n-1}_\infty$).

Thus we are in a situation described in the previous section. For a point
$x\in\CC\PP^{n-1}_\infty\subset\CC\PP^n$, let $p_{\sigma, x}$ be the germ of a 
holomorphic
function (section) $P_\sigma$ at this point (in fact $p_{\sigma, x}$ is a
polynomial in an affine chart there) let ${\hat p}_{\sigma,x}={p_{\sigma,x}}_
{\vert(\CC\PP^{n-1}_\infty, x)}:(\CC\PP^{n-1}_\infty, x)\to\CC$, and let
$\hat{\zeta}_{p_\sigma,x}(t):= {\zeta}_{p_\sigma,x}(t)/{\zeta}_{{\hat p}_\sigma,x}(t)$.
Also let $\hat{\chi}(\{p_{\sigma,x}=0\}):= {\chi}(\{p_{\sigma,x}=0\})
-{\chi}(\{{\hat p}_{\sigma,x}=0\})$. Let ${\widetilde P}_\sigma(x_0, x_1,
\ldots, x_n)$ be the homogenization of degree $d$ of the polynomial
$P_\sigma(x_1, \ldots, x_n)$:
$$
{\widetilde P}_\sigma(x_0, x_1, \ldots, x_n)=x_0^d P_\sigma(\frac{x_1}{x_0},
\ldots, \frac{x_n}{x_0})
$$
(if $\deg P_0 < d$, the polynomial $\widetilde{P}_0$ is not the usual homogenization
of the polynomial $P_0$, but differs from it by the factor $x_0^{d-\deg P_0}$).
Theorem \ref{theo1} gives the following.

\begin{theorem}\label{theo2}
One has
$$
\zeta_{P_\sigma}(t) =
\int_{\{x\in \CC^n:P_0(x)=0\}}\zeta_{P_{\sigma,x}}(t)d\chi\ \cdot
\int_{\{x\in \CC\PP^{n-1}_\infty: {\widetilde P}_0(x)=0\}}
\hat\zeta_{p_{\sigma, x}}(t) d\chi
$$
and therefore
$$
\chi(V_\sigma) = \int_{\{x\in \CC^n: P_0(x)=0\}}\chi(V_{\sigma, x})d\chi\
+ \int_{\{x\in \CC\PP^{n-1}_\infty: {\widetilde P}_0(x)=0\}}
\hat\chi(\{p_{\sigma, x}=0\}) d\chi.
$$
\end{theorem}

\begin{remarks}
{\bf 1.}
If $\deg P_0=d_0<d$, then $\{x\in \CC\PP^{n-1}_\infty: {\widetilde
P}_0(x)=0\}=\CC\PP^{n-1}_\infty$ and at a generic point of the infinite
hyperplane $\CC\PP^{n-1}_\infty$ one has $\widetilde\zeta_{p_{\sigma,
x}}(t)=1-t^{d-d_0}$,
$\widetilde\chi(\{p_{\sigma, x}=0\}_x = d-d_0$.

\noindent{\bf 2.} Theorem~\ref{theo2} reduces computation of the zeta function
of a family of polynomials to computation of the zeta functions of families
of holomorphic germs. An interesting case is a linear family of polynomials
and respectively linear families of holomorphic germs. In somewhat other
terms this case was treated in the study of meromorphic germs
elaborated for study of polynomial maps: \cite{GLM2, GLM3}.
Let $s_\sigma=f+\sigma g$ be a linear family of holomorphic germs (the
family of the zero level sets of $s_\sigma$ is a pencil). Then, modulo the
indeterminacy locus $\{f=g=0\}$, the general local level set and the monodromy
transformation of of the family $s_\sigma$ coincides with the zero Milnor fibre
${\calM}_\varphi^0$ and the corresponding monodromy transformation $h_\varphi^0$
of the meromorphic germ $\varphi=\frac{f}{g}$ as they are defined in
\cite{GLM1}. If $g(0)\ne 0$, the indeterminacy locus is empty and the indicated
objects coincide (and coincide with the usual Milnor fibre of the germ $f$
and its classical monodromy transformation). If $f(0)=g(0)=0$, the indeterminacy
locus is (locally) contractible and the monodromy transformation may be
supposed to be identity on it. Therefore
$\chi(V_{\sigma,0})=\chi({\calM}_\varphi^0)+1$,
$\zeta_{s_\sigma, 0}(t)=\zeta^0_\varphi(t)(1-t)$. 
This permits to apply methods elaborated for meromorphic germs to linear
families of holomorphic germs. In particular, there is a Varchenko type
formula which expresses the zeta function $\zeta^0_\varphi(t)$ in terms
of the Newton diagrams of $f$ and $g$ (in the case when $\varphi$ is
non-degenerate with respect to the Newton diagrams): \cite{GLM2}.

\noindent{\bf 3.} Let $P_\sigma(x)=f(x)+ \sigma g_d(x)$ where
$f$ is a polynomial of degree $d$ such that the projective closure
of any fibre of $f$ and its intersection with the hyperplane at infinity 
have isolated singularities (in fact
isolated boundary singularities) and $g$ is a sufficiently general 
homogeneous polynomial of degree $d$,
such that the compactified fibres of $f+ \sigma g_d$ have transversal 
intersections with infinity as soon as $\sigma \ne 0$.
(e.g. $g= \ell ^d$ for a generic linear function $\ell$).
In \cite[section 7]{ST2} it is shown that the zeta-function of these
deformation is equal to
$$
\zeta_{P_\sigma}(t) = (1-t)^{\chi(V_0)} \prod \zeta_{Z_i}(t)
$$
where $V_0$ is the fibre $\{ f=0 \}$ (supposed to be smooth) and 
$\zeta_{Z_i}$ is the zeta function of the Milnor monodromy of the boundary 
singularities mentioned above.
This formula is now an immediate corrolary of Theorem~\ref{theo2}.
\end{remarks}

\section{Examples}
\noindent{\bf 1.} Let $P_\sigma(x, y)=x^{d_0}+\sigma(x^d+y^d)$ (n=2). There are
3 different cases.

\noindent 1) $d_0>d$. In this case the set $\{\widetilde P_0=0\}\subset\CC\PP^2$
is the closure of the line $\{x=0\}$. There are 3 types of points in it:\newline
a) The origin $0=(0, 0)$. The Varchenko type formula gives
$$\zeta_{P_{\sigma,0}}(t)= (1-t)(1-t^{d_0-d})^{1-d}.$$\newline
b) Other points of the affine line $\{x=0\}$, i.e. $y\ne 0$. At such a point one
has $\zeta_{P_{\sigma,x}}(t)=(1-t^{d_0})$. However, the Euler characteristic of
the set of these points is equal to zero and therefore this stratum gives no
impact to $\zeta_{P_\sigma}(t)$.\newline
c) The infinite point of the line $\{x=0\}$. One can easily see that
$\widetilde\zeta_{p_{\sigma,x}}(t)=1$.

Integrating these local data one gets
$$
\zeta_{P_\sigma}(t)= (1-t)(1-t^{d_0-d})^{1-d}.
$$
One can say that this zeta function essentially originates from the origin in
$\CC^2$.

\noindent 2) $d_0=d$. This case is not interesting : simple computations (or
considerations) give $\zeta_{P_\sigma}(t)= (1-t)$.

\noindent 3) $d_0<d$. In this case the set $\{\widetilde P_0=0\}\subset\CC\PP^2$
is the union of the line $\{x=0\}$ in the affine plane and the infinite line
$\CC\PP^1_\infty$. There are 5 types of points in it:\newline
a) The origin $0=(0, 0)$. One has (e.g., from the Varchenko type formula)
$\zeta_{P_{\sigma,0}}(t)= (1-t)$.\newline
b) Other points of the affine line $\{x=0\}$. Again
$\zeta_{P_{\sigma,x}}(t)=(1-t^{d_0})$, but the Euler characteristic of this
stratum is equal to zero.\newline
c) The infinite point of the line $\{x=0\}$:
$\widetilde\zeta_{p_{\sigma,x}}(t)=1$.\newline
d) Intersection points of the infinite line $\CC\PP^1_\infty$ with the closure 
of the curve $\{x^d+y^d=0\}$; there are $d$ of them. One can easily see that
$\widetilde\zeta_{p_{\sigma,x}}(t)=1$ at them.\newline
e) Finally we have all other (generic) points of the infinite line. The Euler
characteristic of this stratum is equal to $2-(1+d)=1-d$. As it was explaned in
the Remark at the end of Section~\ref{pol},
$\widetilde\zeta_{p_{\sigma,x}}(t)=1-d^{d-d_0}$.

Integrating these local data one gets
$$
\zeta_{P_\sigma}(t)= (1-t)(1-t^{d-d_0})^{1-d}
$$
(almost as in the case 1).
One can say that this zeta function essentially originates from the open stratum
of the infinite line.

\begin{figure}
$$
%% \unitlength=1.00mm
\unitlength=0.50mm
\begin{picture}(120.00,60.00)(-10,-15)
\thinlines
\put(30,30){\line(1,0){39}}
\put(71,30){\vector(1,0){19}}
\put(30,30){\vector(0,1){40}}
\put(30,30){\vector(-1,-1){40}}
\put(30,40){\circle*{1.0}}
\put(30,50){\circle*{2.0}}
\put(30,60){\circle*{1.0}}
\put(40,30){\circle*{1.0}}
\put(50,30){\circle*{1.0}}
\put(60,30){\circle*{1.0}}
\put(70,30){\circle*{1.0}}
\put(80,30){\circle*{1.0}}
\put(90,30){\circle*{1.0}}
\put(60,40){\circle*{2.0}}
\put(70,30){\circle{2.0}}
\put(25,25){\circle*{1.0}}
\put(20,20){\circle*{1.0}}
\put(15,15){\circle*{1.0}}
\put(10,10){\circle*{2.0}}
\put(5,5){\circle*{1.0}}
\thicklines
\put(30,50){\line(3,-1){30}}
\put(30,50){\line(-1,-2){20}}
\put(10,10){\line(5,3){50}}
\put(96,33){$k_2$}
\put(-7,-13){$k_1$}
\put(32,68){$k_0$}
\put(-20,50){{\scriptsize $\zeta_3(t)=(1-t^4)^2$}}
\put(-20,40){{\scriptsize $\zeta_2(t)=1-t^2$}}
\put(-20,30){{\scriptsize $\zeta_1(t)=1$}}
\end{picture}
$$
\caption{The Newton diagrams of the summands corresponding to deformation at the point $(0:0:0:1)$
(marked by $\bullet$ and $\circ$ respectively).}
\label{fig1}
\end{figure}
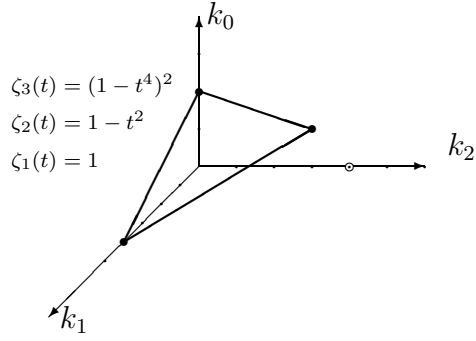

\noindent{\bf 2.} Let $P_\sigma(x_1, x_2, x_3)=x_1^4 + x_2^3 + x_3^2
+ \sigma x_2^4$ (in this example we observe some speical changes at the behaviour at
infinity). The projective variety $\{\widetilde P_0=0\}\subset\CC\PP^3$ intersects
the infinite plane $\CC\PP^2_\infty=\{x_0=0\}$ along the line $\{x_1=0\}$ and
consists of the following (smooth) strata:
\newline
a) The origin $\{0\}$ in $\CC^3$. Here we have an equisingular deformation
of the (surface) singularity $E_6$ and therefore $\zeta_{P_\sigma, 0}(t) = 1-t$.
\newline
b) The set $\{P_0=0\}\setminus\{0\}$. One has $\zeta_{P_\sigma, x}(t) =
1-t$ at all points $x$ of this stratum, however the Euler characteristic
of the stratum itself is equal to zero and thus it does not contribute
to the zeta function of the family $P_\sigma$.
\newline
c) The "distinguished point" $(0:0:0:1)$ on the infinite projective line
$\{x_0=x_1=0\}$. With the help of the Varchenko type formula (\cite{GLM2}; see
Fig.~\ref{fig1})
one gets
$\zeta_{P_\sigma, x}(t) = \frac{(1-t^4)^2}{1-t^2}$ for this point.
\newline
d) The affine line $\{x_0=x_1=0\}\setminus\{(0:0:0:1)\}$. Its Euler
characteristic is equal to $1$. The variety
$\widetilde P_0(0)=0$ is non-singular at points of it, however, its
intersection with the infinite plane consists of a line of multiplicity
$4$ which, for $\sigma\ne0$, splits into $4$ different lines intersecting
each other at one point $(0:0:0:1)$. Therefore, for points of this
stratum, one has $\zeta_{P_\sigma, x}(t) = \frac{1-t}{1-t^4}$.

\noindent
Combining all these local data one gets $\zeta_{P_\sigma}(t)=
\frac{(1-t^4)(1-t)^2}{1-t^2}$.

\noindent{\bf 3.} Let $P_\sigma(x)=f_{d_0}(x)+ \sigma g_d(x)$ where 
$x=(x_1, ..., x_n)\in\CC^n$, $f_{d_0}$ is a non-degenerate homogeneous polynomial of degree $d_0$
(i.e., it has an isolated critical point at the origin),
$g_d$ is a generic homogeneous
polynomial of degree $d$, i.e., the polynomial $g_d$ is non-degenerate and the 
hypersurfaces $\{f_{d_0=0}\}$ and $\{g_d=0\}$ in $\CC\PP^{n-1}_\infty$ intersect
transversally.

\noindent 1) $d_0>d$. There are 3 types of points in the set 
$\{\widetilde P_0=0\}\subset\CC\PP^n$:
\newline a) The origin in $\CC^n$. The Varchenko type formula gives
$$
\zeta_{P_\sigma, 0}(t)=
(1-t)\cdot\left(1-t^{d_0-d}\right)^{(-1)^{n-1}((d_0-1)^n-(d-1)^n)/(d_0-d)}.
$$
\newline b) Other points of the hypersurface $\{P_0=0\}$ in the affine space $\CC^n$.
For these point $\zeta_{P_\sigma, x}(t)=(1-t)$. The Euler characteristic of the set
of these points is equal to zero.
\newline c) Infinite points of the set $\{\widetilde P_0=0\}$, i.e., points of 
$\{f_{d_0=0}\}\subset \CC\PP^{n-1}_\infty$. For these points one has
$\hat\zeta_{p_\sigma, x}(t)=1$.

\noindent
Combining the local data one gets
$$
\zeta_{P_\sigma}(t)=
\zeta_{P_\sigma, 0}(t)=(1-t)\cdot\left(1-t^{d_0-d}\right)^{(-1)^{n-1}((d_0-1)^n-(d-1)^n)/(d_0-d)}.
$$

\noindent 2) $d_0=d$. Obviously $\zeta_{P_\sigma}(t)=1-t$.

\noindent 3) $d_0<d$. The set $\{\widetilde P_0=0\}$ is the union
$\{x\in\CC^n:P_0=0\}\cup\CC\PP^{n-1}_\infty$.
\newline a) For all points of the set $\{x\in\CC^n:P_0=0\}$ one has
$\zeta_{P_\sigma, x}(t)=(1-t)$. The Euler characteristic of this set is equal to 1.
\newline b) At the points of $\{f_{d_0}=0\}\subset\CC\PP^{n-1}_\infty$ and at the points of
$\{g_{d}=0\}\subset\CC\PP^{n-1}_\infty$ one has $\hat\zeta_{p_\sigma, x}(t)=1$.
\newline c) At the points of the complement $\CC\PP^{n-1}_\infty\setminus
(\{f_{d_0}=0\}\cup \{g_{d}=0\})$ one has $\hat\zeta_{p_\sigma, x}(t)=1-t^{d-d_0}$.
The Euler characteristic of this set is equal to $(-1)^{n-1}((d-1)^n-(d_0-1)^n)/(d-d_0)$.

Combying the local data one gets 
$$
\zeta_{P_\sigma}(t)=
(1-t)\cdot\left(1-t^{d-d_0}\right)^{(-1)^{n-1}((d-1)^n-(d_0-1)^n)/(d-d_0)}.
$$

Though for $d_0>d$ and $d_0<d$ the answers are similar, as in Example 1 one can
say that the origins of the zeta function (and/or of the vanishing cycles) in these
cases are different: the origin in $\CC^n$ and the infinite hyperplane $\CC\PP^{n-1}_\infty$ respectively.

\noindent{\bf 4.} Let $P_\sigma(x)=f_{d_0}(x)+ \sigma (\ell(x))^d$ where 
$x=(x_1, ..., x_n)\in\CC^n$, $f_{d_0}$ is a non-degenerate homogeneous polynomial
of degree $d_0$, $\ell$ is a generic (homogeneous) linear function, i.e., the 
hypersurfaces $\{f_{d_0=0}\}$ and $\{\ell=0\}$ in $\CC\PP^{n-1}_\infty$ intersect
transversally. Considerations similar to those of Example 3 give:
$$
\zeta_{P_\sigma}(t)=
\begin{cases}
(1-t)\cdot
(1-t^{d_0-d})^{(1-d_0)^{n-1}} &\mbox{for } d_0 > d,\\
1-t &\mbox{for } d_0=d, \\
(1-t)\cdot
(1-t^{d-d_0})^{(1-d_0)^{n-1}} &\mbox{for } d>d_0.
\end{cases}
$$

\bigskip
\noindent Moscow State University, Faculty of Mechanics and Mathematics\\
Moscow, 119992, Russia\\
E-mail: sabir@mccme.ru

\medskip
\noindent Universiteit Utrechts, Mathematisch Instituut\\
P.O.Box 80.010, 3508 TA Utrecht, The Netherlands\\
E-mail: siersma@math.uu.nl

\end{document}